\documentclass[10pt,a4paper]{article}
\usepackage{latexsym}
\usepackage{amsmath}
\usepackage{graphicx}
\usepackage{amsfonts}
\usepackage{amssymb}
\usepackage{amscd}

\begin{document}
\newtheorem{theorem}{Theorem}[section]
\newtheorem{remark}[theorem]{Remark}
\newtheorem{mtheorem}[theorem]{Main Theorem}
\newtheorem{bbtheo}[theorem]{The Strong Black Box}
\newtheorem{observation}[theorem]{Observation}
\newtheorem{proposition}[theorem]{Proposition}
\newtheorem{lemma}[theorem]{Lemma}
\newtheorem{testlemma}[theorem]{Test Lemma}
\newtheorem{mlemma}[theorem]{Main Lemma}
\newtheorem{note}[theorem]{}
\newtheorem{steplemma}[theorem]{Step Lemma}
\newtheorem{corollary}[theorem]{Corollary}
\newtheorem{notation}[theorem]{Notation}
\newtheorem{example}[theorem]{Example}
\newtheorem{definition}[theorem]{Definition}

\renewcommand{\labelenumi}{(\roman{enumi})}
\newcommand{\dach}[1]{\hat{\vphantom{#1}}}
\newcommand{\cp}{\widehat}
\newcommand{\dsum}{\bigoplus}

\newcommand{\pure}{\subseteq_\ast}
\numberwithin{equation}{section}

\def\Pf{\smallskip\goodbreak{\sl Proof. }}
\def\Fp{\vadjust{}\penalty200 \hfill
\lower.3333ex\hbox{\vbox{\hrule\hbox{\vrule\phantom{\vrule height
6.83333pt depth 1.94444pt width 8.77777pt}\vrule}\hrule}}
\ifmmode\let\next\relax\else\let\next\par\fi \next}

\def\Fin{\mathop{\rm Fin}\nolimits}
\def\br{\mathop{\rm br}\nolimits}
\def\fin{\mathop{\rm fin}\nolimits}
\def\Ann{\mathop{\rm Ann}\nolimits}
\def\Aut{\mathop{\rm Aut}\nolimits}
\def\End{\mathop{\rm End}\nolimits}
\def\bfb{\mathop{\rm\bf b}\nolimits}
\def\bfi{\mathop{\rm\bf i}\nolimits}
\def\bfj{\mathop{\rm\bf j}\nolimits}
\def\df{{\rm df}}
\def\bfk{\mathop{\rm\bf k}\nolimits}
\def\bEnd{\mathop{\rm\bf End}\nolimits}
\def\iso{\mathop{\rm Iso}\nolimits}
\def\id{\mathop{\rm id}\nolimits}
\def\Ext{\mathop{\rm Ext}\nolimits}
\def\Ines{\mathop{\rm Ines}\nolimits}
\def\Hom{\mathop{\rm Hom}\nolimits}
\def\bHom{\mathop{\rm\bf Hom}\nolimits}
\def\Rk{ R_\k-\mathop{\bf Mod}}
\def\Rn{ R_n-\mathop{\bf Mod}}
\def\map{\mathop{\rm map}\nolimits}
\def\cf{\mathop{\rm cf}\nolimits}
\def\top{\mathop{\rm top}\nolimits}
\def\Ker{\mathop{\rm Ker}\nolimits}
\def\Bext{\mathop{\rm Bext}\nolimits}
\def\Br{\mathop{\rm Br}\nolimits}
\def\dom{\mathop{\rm Dom}\nolimits}
\def\min{\mathop{\rm min}\nolimits}
\def\im{\mathop{\rm Im}\nolimits}
\def\max{\mathop{\rm max}\nolimits}
\def\rk{\mathop{\rm rk}}
\def\Diam{\diamondsuit}
\def\Z{{\mathbb Z}}
\def\Q{{\mathbb Q}}
\def\N{{\mathbb N}}
\def\bQ{{\bf Q}}
\def\bF{{\bf F}}
\def\bX{{\bf X}}
\def\bY{{\bf Y}}
\def\bHom{{\bf Hom}}
\def\bEnd{{\bf End}}
\def\bS{{\mathbb S}}
\def\AA{{\cal A}}
\def\BB{{\cal B}}
\def\CC{{\cal C}}
\def\DD{{\cal D}}
\def\TT{{\cal T}}
\def\FF{{\cal F}}
\def\GG{{\cal G}}
\def\PP{{\cal P}}
\def\SS{{\cal S}}
\def\XX{{\cal X}}
\def\YY{{\cal Y}}
\def\fS{{\mathfrak S}}
\def\fH{{\mathfrak H}}
\def\fU{{\mathfrak U}}
\def\fW{{\mathfrak W}}
\def\fK{{\mathfrak K}}
\def\PT{{\mathfrak{PT}}}
\def\T{{\mathfrak{T}}}
\def\fX{{\mathfrak X}}
\def\fP{{\mathfrak P}}
\def\X{{\mathfrak X}}
\def\Y{{\mathfrak Y}}
\def\F{{\mathfrak F}}
\def\C{{\mathfrak C}}
\def\B{{\mathfrak B}}
\def\J{{\mathfrak J}}
\def\fN{{\mathfrak N}}
\def\fM{{\mathfrak M}}
\def\Fk{{\F_\k}}
\def\bar{\overline }
\def\Bbar{\bar B}
\def\Cbar{\bar C}
\def\Pbar{\bar P}
\def\etabar{\bar \eta}
\def\Tbar{\bar T}
\def\fbar{\bar f}
\def\nubar{\bar \nu}
\def\rhobar{\bar \rho}
\def\Abar{\bar A}
\def\a{\alpha}
\def\b{\beta}
\def\g{\gamma}
\def\w{\omega}
\def\e{\varepsilon}
\def\o{\omega}
\def\va{\varphi}
\def\k{\kappa}
\def\m{\mu}
\def\n{\nu}
\def\r{\rho}
\def\f{\phi}
\def\hv{\widehat\v}
\def\hF{\widehat F}
\def\v{\varphi}
\def\s{\sigma}
\def\l{\lambda}
\def\lo{\lambda^{\aln}}
\def\d{\delta}
\def\z{\zeta}
\def\ale{\aleph_1}
\def\aln{\aleph_0}
\def\Cont{2^{\aln}}
\def\nld{{}^{ n \downarrow }\l}
\def\n+1d{{}^{ n+1 \downarrow }\l}
\def\hsupp#1{[[\,#1\,]]}
\def\size#1{\left|\,#1\,\right|}
\def\Binfhat{\widehat {B_{\infty}}}
\def\Zhat{\widehat \Z}
\def\Mhat{\widehat M}
\def\Rhat{\widehat R}
\def\Phat{\widehat P}
\def\Fhat{\widehat F}
\def\fhat{\widehat f}
\def\Ahat{\widehat A}
\def\Chat{\widehat C}
\def\Ghat{\widehat G}
\def\Bhat{\widehat B}
\def\Btilde{\widetilde B}
\def\Ftilde{\widetilde F}
\def\restr{\mathop{\upharpoonright}}
\def\to{\rightarrow}
\def\arr{\longrightarrow}

\newcommand{\norm}[1]{\text{$\parallel\! #1 \!\parallel$}}
\newcommand{\supp}[1]{\text{$\left[ \, #1\, \right]$}}
\def\set#1{\left\{\,#1\,\right\}}
\newcommand{\mb}{\mathbf}
\newcommand{\wt}{\widetilde}
\newcommand{\card}[1]{\mbox{$\left| #1 \right|$}}
\newcommand{\union}{\bigcup}
\newcommand{\inters}{\bigcap}
\def\Proof{{\sl Proof.}\quad}
\def\fine{\ \black\vskip.4truecm}
\def\black{\ {\hbox{\vrule width 4pt height 4pt depth
0pt}}}
\def\fine{\ \black\vskip.4truecm}
\long\def\alert#1{\smallskip\line{\hskip\parindent\vrule%
\vbox{\advance\hsize-2\parindent\hrule\smallskip\parindent.4\parindent%
\narrower\noindent#1\smallskip\hrule}\vrule\hfill}\smallskip}

\title{\sc Absolutely Indecomposable Modules}

\footnotetext{This work is supported by the project No.
  I-706-54.6/2001 of the German-Israeli
  Foundation for Scientific Research \& Development.\\
  This is GbSh880 in second author's list of publications.
  }

\author{R\"udiger G\"obel and Saharon Shelah}

\date{}

\maketitle

\begin{abstract} A module is called {\sl absolutely indecomposable}
if it is directly indecomposable in every generic extension of the
universe. We want to show the existence of large abelian groups that
are absolutely indecomposable. This will follow from a more general
result about $R$-modules over a large class of commutative rings $R$
with endomorphism ring $R$ which remains the same when passing to a
generic extension of the universe. It turns out that `large' in this
context has a {\em precise meaning,} namely being smaller than the
first $\o$-Erd\H{o}s cardinal defined below. We will first apply a
result on large rigid valuated trees with a similar property
established by Shelah \cite{S} in 1982, and will prove the existence
of related `$R_\o$-modules' ($R$-modules with countably many
distinguished submodules) and finally pass to $R$-modules. The
passage through $R_\o$-modules has the great advantage that the
proofs become very transparent essentially using a few `linear
algebra' arguments accessible also for graduate students. The result
closes a gap in \cite{ES,EM}, provides a good starting point for
\cite{FG} and gives a new construction of indecomposable modules in
general using a counting argument.
\end{abstract}

\section{Introduction}

There is a whole industry transporting symmetry properties from one
category to another: For example consider a tree or a graph (with
extra properties if needed) together with its group of
automorphisms. Then encode the tree or the graph into an object of
your favored category in such a way that the branches (or vertices)
of the tree (of the graph) are recognized in the new structure. If
the new category are abelian group argue by (infinite) divisibility,
in case of groups and fields you use of course infinite chains of
roots (with legal primes) etc. Thus the automorphism group of the
tree or the graph is respected in the new category and by density
arguments (or killing unwanted automorphisms by prediction arguments
`on the way') it happens that the automorphism group we start with
becomes (modulo inessential maps: inner automorphisms in case of
groups and Frobenius automorphisms in case of fields) the
automorphism group of an object of the new category. For a few
illustrating details the reader may want to see papers by Heineken
\cite{H}, Braun, G\"obel \cite{BG1} (in case of groups), Corner,
G\"obel \cite{CG} in case of modules (with group rings as the first
category) or Fried and Kollar \cite{FK}, Dugas, G\"obel \cite{DG} in
case of fields and \cite{DG1} for automorphism groups of geometric
lattices. In this paper we also argue with symmetry properties of
trees, but they are of a different kind. Given a cardinal $\l$ which
is not extraordinary large (we explain what we mean by
`extraordinary large' in the next section) then there is an absolute
and rigid family of ($\o$-)valuated trees based on this cardinal.
This is a family of $\l$ subtrees $T$ of size $\l$ of the tree
$T_\l= {}^{\o>}\l$ of finite sequences of ordinals in $\l$ together
with a valuation map $v:T\arr \o$. Rigid means that there is no
level preserving valuated homomorphisms between any two distinct
members. (A tree homomorphism is valuated if the value of a branch
is the same as the value of its image.) Moreover this property is
preserved if we change the universe, passing to a generic extension
of the given universe (of set theory) we live in. The existence of
such trees was shown by Shelah \cite{S}. These trees (used also in
applied mathematics) were considered earlier in papers by
Nash-Williams, see \cite{Na} for example. We will encode them into
free $R_\o$-modules over an arbitrary not extraordinary large
commutative ring $R$ with $1\ne 0$. To be definite we can assume
that $R$ is the field $\Q$ of rationals or $\Z$. Recall that
$R_\o$-modules are $R$-modules with countably many ($\o$)
distinguished submodules and free means that the module and its
distinguished submodules and factor modules are free as well. Such
creatures are considered in Brenner, Butler, Corner (see
\cite{B,BB,C,BG}) and G\"obel, May \cite{GM} for arbitrary
commutative rings and an account about the advanced theory in case
of fields can be seen in \cite{Si} and in the references given
there. We will show the existence of free $R_\o$-modules with
endomorphism algebra $R$ by transporting the absolute rigid trees
into the category of $R_\o$-modules. It turns out that the passage
through $R_\o$-modules makes the anticipated proofs very
transparent. Moreover our main result on $R_\o$-modules with
distinguished submodules is only a few steps away from the desired
result on $R$-modules if $R$ has enough primes (like $\Z$).

The corollary on the existence of large absolutely (fully) rigid
abelian groups replaces the earlier unsuccessful approach in
\cite{ES} and \cite[Chapter XV]{EM}: Let $R\ne 0$ be any fixed
countable ring. Then by Corollary \ref{kappaomega} there exists an
absolutely rigid $R_\o$-module of size $\l$ (or an absolute family
of size $\l$ of non-trivial $R$-modules with only the
zero-homomorphism between distinct member) iff $\l < \k(\o)$. The
same holds if $R_\o$-modules are replaced by abelian groups. Thus as
a byproduct we present a new construction of large, absolutely
indecomposable abelian groups, not using stationary sets as
\cite{S0,CG}. So, if we restrict to the problem on the existence of
large absolute indecomposable abelian groups addressed in
\cite{ES,EM}, then it follows from the above (realizing for example
$\Z$ as the endomorphism ring in Corollary \ref{kappaomega}) that
from $\l < \k(\o)$ follows the existence of such abelian groups. The
converse direction would need a strengthening of the Theorem
\ref{nonscalar} from \cite{ES} now showing the existence of
non-trivial idempotents. (The second author believes that this guess
might be true.)

It is also a different matter how to replace $R_\o$-modules by
$R_4$-modules or $R_5$-modules and the endomorphism algebra $R$ by a
general not extraordinary large prescribed $R$-algebra $A$. This
will follow from \cite{FG}, a paper which had to wait for Theorem
\ref{fully} in place of \cite{ES}.

\bigskip

\section {Rigid families of valuated trees and the first $\boldsymbol \o$--Erd\H{o}s
cardinal}

We first describe the result on trees we want to apply by encoding
them into modules with distinguished submodules.

Let $\k(\o)$ denote the first $\o$-{\sl Erd\H{o}s cardinal}. This is
defined as the  smallest cardinal $\k$ such that $\k \to (\o)^
{<\o}$, i.e. for every function $f$  from the finite subsets of $\k$
to 2 there exist an infinite subset $X \subset \k$ and a function
$g: \o \to 2$ such that $f(Y) = g(|Y|)$ for all finite subsets $Y$
of $X$. The cardinal $\k(\o)$ is strongly inaccessible; see Jech
\cite[p. 392]{J}. Thus $\k(\o)$ is a large cardinal. We should also
emphasize that $\k(\o)$ may not exist in every universe.

If $\l <\k(\o)$, then let $$T_\l= {}^{\o>}\l=\{f: n\arr \l: \text{
with } n < \o \text{ and } n = \dom f  \}$$ be the tree of all
finite sequences $f$ (of length or level $\lg(f)=n$) in $\l$.
Since $n =\{0,\dots n-1\}$ as ordinal, we also write
$f=f(0){}^\wedge f(1){}^\wedge \dots {}^\wedge f(n-1)$. By
restriction $g=f\restr m$ for any $m\le n$ we obtain all {\em
initial segments} of $f$. We will write $g\lhd f$. Thus
$$g\le f\iff g\subseteq f \text{ as graphs }\iff g\lhd f.$$
A subtree $T$ of $T_\l$ is a subset which is closed under initial
segments and a homomorphism between two subtrees of $T_\l$ is a map
that preserves levels and initial segments. (Note that a
homomorphism does not need to be injective or preserve $\nleq$.) The
tree $T$ is valuated if with the tree we have a valuation map
$v:T\arr \o$. In the following a tree will always come with a
valuation and $\Hom(T_1,T_2)$ denotes the valuated homomorphisms
between subtrees $T_1$ and $T_2$, i.e if $v_i$ is the valuation of
$T_i$ ($i=1,2$) and $\va$ is such a \emph{valuated homomorphism},
then $v_2(\eta\va) = v_1(\eta)$ for all $\eta\in T_1$. Shelah
\cite{S} showed the existence of an {\em absolutely rigid } family
of $2^\l$ valuated subtrees of $T_\l$.
\begin{theorem} \label{sheltrees} If $\l<\k(\o)$ is infinite and
$T_\l= {}^{\o>}\l$, then there is a family $T_\a$ ($\a\in 2^\l$) of
valuated subtrees of $T_\l$ (of size $\l$) such that for $\a,\b\in
2^\l$ and in any generic extension of the universe the following
holds.
$$\Hom(T_\a,T_\b)\ne \emptyset \Longrightarrow \a=\b.$$
\end{theorem}

\Pf The result is a consequence of the Main Theorem 5.3 in \cite[p.
208]{S}. The family of rigid trees is constructed in \cite[p. 214,
Theorem 5.7]{S} and the proof, that the trees are rigid, follows
from Theorem 5.8 using the Conclusion 2.14 in \cite{S}. In Shelah's
notation $\k(\omega)$ is the first beautiful cardinal $>\aleph_0$.
\fine

This property of rigid families of valuated trees in Theorem
\ref{sheltrees} fails, if we choose $\l\ge\k(\o)$. In fact the
following result from \cite{ES} on rigid families of $R$-modules
reflects this.

\begin{theorem}\rm{(Eklof-Shelah \cite{ES})}\label{nonscalar} Let
$\l$ be a cardinal $\geq \k(\o)$ and $R$ any ring with $1$.
\begin{enumerate}
    \item If $\{M_\nu  \ |\ \nu  < \l\}$ is a family of non-zero left
    $R$-modules, then there are distinct ordinals $\mu, \nu < \l$,
    such that in some generic extension $V[G]$ of the universe $V$,
    there is an injective homomorphism  $\phi: M_\mu \to M_\nu$.
    \item If $M$ is an $R$-module of cardinality $\l$, then there exists
    a generic extension  $V[G]$ of the universe $V$, such that $M$
    has an endomorphism that is not multiplication by an element of $R$.
\end{enumerate}
 \end{theorem}

Thus $\k(\o)$ is the precise border line for Theorem \ref{rep} and
we can not expect absolute results on endomorphism rings and rigid
families of abelian groups above $\k(\o)$, see Corollary~
\ref{kappaomega}.

Combining Theorem \ref{nonscalar} with our main result this also
conversely shows that the implication of Theorem \ref{sheltrees}
fails whenever $\l\ge\k(\o)$, i.e. there is a generic extension
$V[G]$ of the universe $V$ and there are distinct ordinals
$\a,\b\in 2^\l$ with $\Hom(T_\a,T_\b)\ne \emptyset$.

\section {The main construction}

Let $R\ne 0 $ be a commutative ring. As we shall write endomorphisms
on the right, it will be convenient to view $R$-modules as left
$R$-modules. Next we define a free $R$-module $F$ of rank $\l$ over
a suitable indexing set (obviously) used to encode trees $T_\a$ from
Theorem~ \ref{sheltrees} into the structure when turning the free
$R$-module $F$ into an $R_\o$-module module $\bF$ with $\o$
distinguished submodules.

We enumerate a subfamily of $\lambda$ valuated trees from Theorem
\ref{sheltrees} by the indexing set $I =  {}^{\omega >}({}^{\omega
> }\lambda)$. Thus
$$ T_{\etabar} \text{ with valuation map } \ v_{\etabar}:T_{\etabar}\arr \o\
 \ (\etabar \in I)$$
without repetition. Next define inductively subsets $S_n \subseteq
{}^n({}^{\omega > }\lambda)$ such that the following holds.
\begin{enumerate}
    \item [(0)] $S_0 = \{\bot\}$
    \item [(1)] If $S_n$ is defined, then $S_{n+1} =
    \{\etabar^\wedge \langle \nu \rangle : \etabar \in S_n, \bot \ne \nu
    \in T_{\etabar} \}$.
\end{enumerate}

Let $S = \bigcup\limits_{n\in \omega}  S_n$ and also let $\etabar
^\wedge \langle \bot \rangle = \etabar$ for $\bot \in T_{\etabar}$.

Put  $S_{nk} = \{\etabar ^\wedge \langle \nu \rangle \in S :\ \lg
\etabar = n, \lg \nu = k \} \subseteq S_{n+1}$. Here $\nu =
\nu_0{}^\wedge \dots {}^\wedge \nu_{k-1}$ with $\nu_i \in \ \l$ is
a sequence of ordinals and $\etabar = \etabar_0{}^\wedge \dots
{}^\wedge \etabar_{n-1}$ with $\etabar_i \in
T_{{\etabar}_0{}^\wedge \dots {}^\wedge \etabar_{i-1}}$ a sequence
of branches from special trees. Moreover write
$$T^k_{\etabar} = \{ \nu \in T_{\etabar} :\ \lg \nu = k \} \subseteq
T_{\etabar}  \text{ and } T_{k{\etabar}} = \{ \nu \in T_{\etabar} :\
v_{\etabar}( \nu) = k \}\ \ \ (\etabar \in I).$$ Now we define the
free $R$-modules:
\begin{enumerate}\label{freemod}
    \item $ F = \bigoplus\limits_{\etabar \in S} R e_{\etabar}$
    \item $F_{nk} = \bigoplus\limits_{\etabar \in S_n}
    \bigoplus\limits_{\nu \in T^k_{\etabar}} R ( e_{\etabar^\wedge
    \langle \nu \restr k-1 \rangle} - e_{\etabar^\wedge \langle \nu \rangle})$
    \item $F^{nk} = \bigoplus\limits_{\etabar \in S_n} \ \ (
    \bigoplus\limits_{\nu \in T^k_{\etabar}}  Re_{\etabar^\wedge \langle \nu
    \rangle})$
    \item $F^k_n = \bigoplus\limits_{\etabar \in S_n} \ \ (
    \bigoplus\limits_{\nu \in T_{k{\etabar}}}  Re_{\etabar^\wedge \langle \nu
    \rangle})$
    \item $F_0 = \langle R(e_{\etabar}\ - e_{\etabar'}) : \
    \etabar, \etabar' \in S \rangle$ and $ F_1 = R e_\bot$.
\end{enumerate}

We note that $F_0 = \bigoplus\limits_{\bot \, \ne \, \etabar\,
\in\, S}R(e_\bot  - e_{\etabar})$ and $F = F_0 \oplus F_1$.

Next we define $R_\o$-modules. These are $R$-modules with $\o$
distinguished submodules. We enumerate the distinguished
submodules by a particular well-ordered, countable indexing set
$$ W= \langle 0,1\rangle^\wedge L_1{}^\wedge L_2{}^\wedge L_3 \text{ with } L_i \text{ a copy of }
\o\times \o \ \ \ (i=1,2,3).$$

We view $W$ as an ordinal. Then an {\em $R_\o$-module $\bX$} is an
$R$-module $X$ with a family of submodules $X_i \ (i\in W)$. We will
also say that $\bX$ is {\em a free $R_\o$-module} if $X,X_i,X/X_i \
(i\in W)$ are free $R$-modules. In particular

 \begin{eqnarray}\label{theF} \bF = (F, F_0,F_1, F_{nk}, F^{pq}, F^r_s: \
 (nk)\in L_1, (pq)\in L_2), (rs)\in L_3) \text{ is a free }
 R_\o-\text{module.}\end{eqnarray}

 If $\bX,\bY$ are $R_\o$-modules, then $\va$ is an $R_\o$-homomorphism
 ($\va\in \bHom_R(\bX,\bY)$) if $\va \in \Hom_R(X,Y)$ and $X_i\va\subseteq Y_i$ for all
 $i\in W$, where $\bY = (Y, Y_i : \ i\in W)$. We also write
 $\bHom_R(\bX,\bX)=\bEnd_R\bX$.

We want to show the following
 \begin{theorem} \label{rep} Let $R$ be a
commutative ring with $1 \neq 0$ and $|R|, \l < \k (\o)$. A free
$R$-module $F$ of rank $\l$ can be made into a free $R_\o$-module
$\bF = (F, F_i :\ i \in W)$ such that $\bEnd_R \bF = R$ holds in
any generic extension of the given universe.
 \end{theorem}

Note that the size of $R$ and the rank $\l$ can be arbitrary $<\k
(\o)$; in particular $R=\Z/2\Z$. If $\l$ is finite, then we can
choose directly a suitable finite family of $F_i$s with the required
endomorphism ring. Otherwise $\l$ is infinite and we can apply
Theorem \ref{sheltrees}. So we choose $\bF = (F, F_i : i \in W)$ as
in (\ref{theF}) depending on the valuated trees from Theorem
\ref{sheltrees}. Then clearly it remains to show $\bEnd_R \bF = R$.
We first show the following crucial

\begin{lemma}\label{diag} Let $\va \in \bEnd_R \bF$ with $\bF$ as in (\ref{theF})
and $F=\bigoplus\limits_{\etabar \in S} R e_{\etabar} $. If
$\etabar\in S$, then $$e_{\etabar}\va \in Re_{\etabar}.$$
\end{lemma}
\Pf  Let $\etabar \in S$ be fixed and recall that $T^k_{\etabar} =
T_{\etabar} \cap {}^k\l$. We consider its successors $\etabar
^\wedge \langle \nu \rangle$ in $S$ with $\bot \neq \nu \in
T_{\etabar}$ and let $\lg \etabar = n, \lg \nu = k$. Thus $\etabar
^\wedge \langle \nu \rangle \in S_{n k}$ and $\nu \in
T^k_{\etabar}$. If $\varphi \in \bEnd_R \bF$, then we claim

\begin{eqnarray} \label{formel} e_{\etabar{}^\wedge \langle \nu \rangle} \va =
\sum\limits_{l < l_\nu}\, \, r_{\nu l} e_{\rhobar_{\nu l}{}^\wedge
\langle \s_{\nu l} \rangle}  \text{ with } \rhobar_{\nu l} \in S_n,
\ \sigma_{\nu l} \in T^k_{\rhobar_{\nu l}} \text{ and } 0 \ne r_{\nu
l} \in R.
\end{eqnarray}

If $e_{\etabar{} ^\wedge \langle \nu \rangle} \varphi = 0$, we
choose $l_\nu=0$ and have the empty sum which is $0$. By definition
of $F^{nk}$ follows $e_{\etabar{}^\wedge \langle \nu \rangle} \in
F^{nk}$, thus $e_{\etabar^\wedge \langle \nu \rangle} \varphi \in
F^{nk}$ showing that $e_{\etabar^\wedge \langle \nu \rangle}\va $ is
of the desired form (\ref{formel}).

We will now use $\bF$ to derive further restrictions of the
expressions in (\ref{formel}).

If $\nu_1 \in T^{k+1}_{\etabar}$, then $\nu_0 = \nu_1 \restr k \in
T^k_{\etabar}$ and $e_{\etabar ^\wedge \langle \nu_0 \rangle} -
e_{\etabar ^\wedge \langle \nu_1 \rangle} \in F_{n\ k+1}$ hence $w
:= (e_{\etabar ^\wedge \langle \nu_0 \rangle} - e_{\etabar ^\wedge
\langle \nu_1 \rangle})\va \in F_{n\ k+1}$ as well. Using
(\ref{formel}) and the definition of $F_{n\ k+1}$ we get
$$ w = \sum\limits_{l < l_{\nu_0}} \ r_{{\nu_{_0}} l} e_{\rhobar_{ \nu_0 l} {}^\wedge
\langle \sigma_{{\nu_0} l} \rangle}  - \sum\limits_{l < l_{\nu_1}}
r_{{\nu_1}l} e_{\rhobar_{\nu_1 l}{}^\wedge \langle \s _{{\nu_1
l}}\rangle}
 = \sum\limits_{i < l_w}  s_{w i}\left( e_{\rhobar_{w i} {}^
\wedge \langle \nu_{w i} \restr \ k \rangle} - e_{\rhobar_{w i} {}^
\wedge \langle \nu_{w i}\rangle} \right)$$
 with $\rhobar_{w i} \in S_n, \ \nu_{w i} \in T^{k+1}_{\rhobar_{w i}}$ and $0 \ne s_{w i} \in R$.

 Now we collect terms of length $k$ and $k+1$ respectively, and
it follows

 $$\ \ \ \ \ \text{ length k: }  \  \ \sum\limits_{l < l_{\nu_0}} \
 r_{{\nu_{_0}}l} e_{\rhobar_{ \nu_0 l} {}^\wedge
\langle \sigma_{{\nu_0} l} \rangle} = \sum\limits_{i < l_w} s_{w i}
e_{\rhobar_{w i} {}^ \wedge \langle \nu_{w i} \restr \ k \rangle}$$
 $$\text{ length k+1: }  \  \sum\limits_{l < l_{\nu_1}} r_{{\nu_1}l}
 e_{\rhobar_{\nu_1 l}{}^\wedge \langle
\s _{{\nu_1 l}}\rangle} =\sum\limits_{i < l_w} s_{w i} e_{\rhobar_{w
i} {}^ \wedge \langle \nu_{w i}\rangle} .$$

We will apply the two displayed equations and suppose for
contradiction that $e_{\etabar}\va \notin R e_{\etabar} $. Hence
$e_{\etabar}\va =\sum\limits_{l < l_{\etabar}} r_l\ e_{\etabar_l}$
and there is $\etabar_0\ne \etabar$ with $r_0\ne 0$. We want to
construct a (level preserving) valuated homomorphism
$$ g: T_{\etabar}\arr T_{\etabar_0} \text{ with } v_{\etabar_0}(g(\nu))
= v_{\etabar}(\nu) \text{ for all } \nu \in T_{\etabar}.$$ Hence
$T_{\etabar},T_{\etabar_0}$ are not rigid and this would contradict
the implication of Theorem \ref{sheltrees}. We will construct
$g=\bigcup_{k\in\o} g_k$ as the union of an ascending chain of
valuated homomorphisms $$g_k: T_{\etabar}\cap {}^{k \ge}\l \arr
T_{\etabar_0}\cap {}^{k \ge}\l.$$

Let  $g_0(\bot) =\bot$ and suppose that $g_k$ is defined subject
to the following condition which we carry on by induction.

\begin{eqnarray}\label{bag} \text{If } \nu_1\in T^k_{\etabar},
\text{ then } \etabar_0{}^\wedge \langle g_k(\nu_1)) \in
\{\rhobar_{\nu_1l}{}^\wedge \langle \s_{\nu_1l}\rangle :\ l <
l_{\nu_1}\}\end{eqnarray}
 thus $g_k(\nu_1)\in T_{\etabar_0}$ for $\etabar_0=
 \rhobar_{\nu_1l}$. Note that (\ref{bag}) is satisfied for $k=0$
 by the assumption on $\va$. Thus we can proceed. If now
 $\nu_1\in T_{\etabar}^{k +1}$ and $\nu_0=\nu_1\restr k$,
 then $g_k(\nu_0)\in T_{\etabar_0}^k$ is given and we want to
 determine $g_{k+1}(\nu_1)$. By induction hypothesis we have some
 $l_*<l_{\nu_0}$ with $\rhobar_{\nu_0l_*}=\etabar_0$ and
 $g_k(\nu_0)=\s_{\nu_0l_*}\in T_{\etabar_0}$.

We must find $l'< l_{\nu_1}$ (see (\ref{formel})) such that
$\rhobar_{\nu_1l'}=\etabar_0$ and $\s_{\nu_0l_*}=\s_{\nu_1l'}\restr k$.
The second condition ensures that $g$ will be the union of an ascending
chain of $g_k's$  and also level preserving. The first assertion
is our induction-bag which we must carry along. It is also the link to
the undesired map $\va$.

By the displayed equation for length $k$, there is some $i$ (perhaps
more than one) such that $s_{w i} \ne 0$ and $\rhobar_{w_i}
=\etabar_0$ and $\nu_{w i}\restr k = \s_{\nu_0l^*}$. Then the other
displayed equation of length $k+1$, by picking one of the preceding
$i$, yields the desired $l'$.

We now have $l'< l_{\nu_1}$ with $\rhobar_{\nu_1l'}=\etabar_0$ and
$\s_{\nu_1l'}\in T_{\etabar_0}$ of length $k+1$ with
$\s_{\nu_1l'}\restr k =\s_{\nu_0l_*}.$ So we can map
$g_{k+1}(\nu_1)\in T_{\etabar_0}$. If $v_{\etabar}(\nu_1) = k$, then
(using $\lg(\etabar)=n$)  $e_{\etabar ^\wedge \langle \nu_1
\rangle}\in F_n^k$ and by (iv) also $e_{\etabar ^\wedge \langle
\nu_1 \rangle}\va \in F_n^k$ and $v{\etabar_0}(g_{k+1}(\nu_1)) = k =
v{\etabar}(\nu_1)$ follows. Thus valuation is preserved.

We argue like this for all $\nu_1\in T_{\etabar}$ of length $k+1$.
This completes the definition of $g_{k+1}$. Thus $g:T_{\etabar}\arr
T_{\etabar_0}$ exists, a contradiction. \fine

\Pf {\sl(of Theorem \ref{rep})} From Lemma \ref{diag} follows
$e_\bot\va = r e_\bot,\ e_{\etabar}\va =r_{\etabar} e_{\etabar} $
for some $r,r_{\etabar}\in R$ and all $\bot\ne \etabar\in S$.
Moreover $(e_\bot -e_{\etabar})\in F_0$, and therefore $(e_\bot
-e_{\etabar})\va \in F_0$ and $(e_\bot -e_{\etabar})\va= r e_\bot
-r_{\etabar} e_{\etabar} \in R(e_\bot -e_{\etabar})$ by support
(in the direct sum). Hence $r e_\bot -r_{\etabar} e_{\etabar}=r'
(e_\bot -e_{\etabar})$ for some $r'\in R$ and $r=r', r_{\etabar}=
r'$ implies $r_{\etabar}= r$ for all $\etabar\in S$. Thus $\va =
r\in R$. \fine

\section{Extension to fully rigid systems}

We want to strengthen Theorem \ref{rep} showing the existence of
fully rigid systems of $R_\o$-modules on $\l$. This is a family
$\bF_U \ (U\subseteq \l)$ of $R_\o$-modules such that the
following holds.

$$
\bHom_R (\bF_U, \bF_V) = \left\{
\begin{array}{rr}
R  \text{  if   }\quad  U \subseteq V\\
0 \text{  if   }\quad  U \not\subseteq V
\end{array} \right.
$$

This result will be the starting point for realizing $R$-algebras
$A$ as endomorphism algebras $\bEnd_R\bF=A$ which are also absolute,
see Fuchs, G\"obel \cite{FG}. We first extend the well-ordered
indexing set $W$ for $\bF$ by one more element and let $$W':=
\langle 0,1,2\rangle ^\wedge L_1{}^\wedge L_2{}^\wedge L_3 \text{
with } L_i\cong \o\times\o.$$
 Hence $W'$ and $W$ are both order-isomorphic to $\o\times\o$ but
$W'$ has virtually one more element than $W$ added at place $2$ to
the definition of $\bF$. This allows us to replace $\bF$ from
Theorem \ref{rep} by $\bF_U$ where the new place is
 $$F_2 := F_U:= \bigoplus\limits_{e\in U} e R \text{ for any }
 U\subseteq S.$$
 From Theorem \ref{rep} follows
 $$\bHom_R(\bF_U.\bF_V)\subseteq R \text{ for any } U,V \subseteq
 S.$$
 Clearly $\bHom_R(\bF_U.\bF_V)= R$ if $U\subseteq V$. On the
 other hand, if $u\in U\setminus V$, then $e_u\va =r e_u $ by the
 displayed formula. But $r e_u \in F_V$ only if $r=0$. Hence
 $\bHom_R(\bF_U.\bF_V)=0$ whenever $U\not\subseteq V$. Finally
 note that $|S|=\l$. We established the existence of fully rigid
 systems.

 \begin{theorem}\label{fully} If $R$ is any commutative ring with
 $1\ne 0$ and $\l,|R|<\k(\o)$, then there is a
 fully rigid system $\bF_U \ (U\subseteq \l)$ of free $R_\o$-modules
 with the following properties.
 \begin{enumerate}
    \item $F$ is free of rank $\l$ and
    $\bF_U=(F,F_0,F_1,F_U, F_i :\ i\in L_1{}^\wedge L_2{}^\wedge L_3)$,
 thus only $F_2=F_U$ depends on
$U$.
    \item The family $\bF_U \ (U\subseteq \l)$ is absolute, i.e.
    if the given universe is replaced by a generic extension, then
    the family is still fully rigid.
 \end{enumerate}
\end{theorem}

The last theorem and a result from \cite{ES} (see Theorem
\ref{nonscalar}) immediately characterize the first $\o$-Erd\H{o}s
cardinal. For clarity we restrict ourself to countable rings $R$.

\begin{corollary}\label{kappaomega} Let $R$ by any countable commutative ring. Then the
following conditions for a cardinal $\l$ are equivalent.
\begin{enumerate}
    \item There is an absolute $R_\o$-module $\bX$ of size $\l$
    with $\End_RM=R$.
    \item There is a fully rigid family $\bF_U \ (U\subseteq \l)$ of
    free $R_\o$-modules.
    \item There is a family of $R_\o$-modules of size $\l$ with
    only the zero-homomorphism between two distinct members.
    \item $\l < \k(\o)$ with $\k(\o)$ the first $\o$- Erd\H{o}s cardinal.
\end{enumerate}
\end{corollary}

We note, that the last theorem can also be applied to vector spaces
(and $\omega$ in (i), (ii) and (iii) can be replaced by $4$ or $5$
as demonstrated in \cite{FG})

\section{ Passing to $R$-modules}
We will restrict ourself to only one application of Theorem
\ref{fully}. A forthcoming paper by Fuchs, G\"obel \cite{FG} will
exploit Theorem \ref{fully} and new results will be obtained in
two directions. Firstly the number of primes needed in Corollary
\ref{infprimes} will be reduced to four (which is minimal),
moreover $R$-algebras $A$ will be realized as $\End_RM=A$ in order
to give more absolute results. These applications were obtained
earlier but had to wait for publication until it became possible
to replace certain results in \cite{ES} by Theorem \ref{fully}.

\begin{corollary} \label{infprimes} Let $R$ be a domain
with infinitely many comaximal primes. If $\l,|R|<\k(\o)$, then
there is an absolute fully rigid family $M_U \ (U\subseteq \l)$ of
torsion-free $R$-modules $M_U$ of size $\l$. Thus the following
holds in any generic extension of the given universe of set
theory.
$$
\Hom_R (M_U, M_V) = \left\{
\begin{array}{rr}
R  \text{  if   }\quad  U \subseteq V\\
0 \text{  if   }\quad  U \not\subseteq V
\end{array} \right.
$$
\end{corollary}

\Pf Let $p_i \ (i\in W')$ be a countable family of comaximal primes
of $R$ and choose $\bF_U= (F,F_0,F_1,F_U, F_i :\ i\in L_1{}^\wedge
L_2{}^\wedge L_3)$ from Theorem \ref{fully}. We will now construct
$R$-modules $M_U$ with
$$F\subseteq M_U\subseteq Q\otimes F$$
where $Q$ denotes the quotient field of $R$. Also, if $X\subseteq
F$, then we denote by
$$p^{-\infty}X:=\bigcup_{n\in\o}p^{-n}X\subseteq Q\otimes F.$$
Now let $$M_U := \langle p_i^{-\infty}F_i, p_2^{-\infty}F_U: i\in
W \rangle.$$

Thus $F\subseteq M_U\subseteq Q\otimes F$ because $F_0+F_1=F$ and
$Q\otimes \bF_U := (Q\otimes F,Q \otimes F_0,Q \otimes F_1, Q\otimes
F_U, Q \otimes F_i :\ i\in L_1{}^\wedge L_2{}^\wedge L_3)$ satisfies
$\End_Q(Q \otimes \bF_U) =Q$ by Theorem \ref{fully}. Consider now
any $\va\in \End_RM_U$. The primes ensure that
$p_i^{-\infty}F_i\va\subseteq p_i^{-\infty}F_i$ for all $i\in W'$
and $\va$ extends uniquely to an endomorphism (also called) $\va\in
\End_Q (Q\otimes \bF_U)$. It follows that $\va=q\in Q$, thus $\va$
is scalar multiplication by $q$ on the right. It remains to show
that $(\va=) q \in R$ and possibly $\va=0$.

Now we recall that the family of primes, in particular $p_0$ and
$p_1$ are comaximal, thus $p_0^{-\infty}R\cap p_1^{-\infty}R = R$.
Choose any $e_{\etabar}\in F_1$. Then $e_{\etabar}\va \in
p_1^{-\infty}R e_{\etabar}$, hence $q\in p_1^{-\infty}R$. Similarly,
$e_\bot\va\in p_0^{-\infty}Re_\bot$, thus also $q\in p_0^{-\infty}R$
and $q\in p_0^{-\infty}R\cap p_1^{-\infty}R = R$ as required. If
$U\not\subseteq V$, then $\Hom_Q(Q\otimes F_U,Q\otimes F_V)=0$ by
Theorem \ref{fully} and the unique extension of $\va$ to the
corresponding $Q$-vector space must by zero. Hence $\va=0$ and the
corollary follows. \fine

We would like to mention that the infinite set of primes in the
corollary can be replaced by $4$ primes, see \cite{FG}; and primes
can also be replaces by comaximal multiplicatively closed subsets.
The latter is a natural straight extension suggested by Tony
Corner (unpublished); this can be looked up in \cite{FrG}.

\noindent Address of the authors:

\noindent
R\"udiger G\"obel \\ Fachbereich 6, Mathematik\\
Universit\"at Duisburg Essen\\ D 45117 Essen, Germany \\
{\small e-mail: r.goebel@uni-essen.de}
\medskip

\noindent Saharon Shelah\\
Institute of Mathematics,
\\Hebrew University, Jerusalem, Israel \\
and Rutgers University, New Brunswick, NJ, U.S.A \\
{\small e-mail: shelah@math.huji.ac.il}
\medskip


\begin{thebibliography}{99}\label{lit} \markright{}


\bibitem{BG} C. B\"ottinger and R. G\"obel, Endomorphism algebras
of modules with distinguished partially ordered submodules over
commutative rings, J. Pure Appl. Algebra {\bf 76} (1991), 121 --
141.

\bibitem{BG1} G. Braun and R. G\"obel, Outer automorphism groups of
locally finite $p$-groups, J. Algebra {\bf 264}, (2003) 55 -- 67.


\bibitem{B} S. Brenner, Endomorphism algebras of vector spaces
with distinguished sets of subspaces, J. Algebra {\bf 6} (1967), 100
-- 114.

\bibitem{BB} S. Brenner and M.C.R. Butler, Endomorphism rings of
vector spaces and torsion free abelian groups, J. London Math. Soc.
{\bf 40} (1965), 183 -- 187.


\bibitem{C} A. L. S. Corner, Endomorphism algebras of large modules
with distinguished submodules, J. Algebra {\bf 11} (1969), 155 --
185.


\bibitem{C1} A. L. S. Corner,  Fully rigid systems of modules, Rendiconti Sem.
 Mat. Univ. Padova {\bf 82} (1989), 55 -- 66.

\bibitem{CG}  A. L. S. Corner and R. G\"obel, Prescribing endomorphism algebras
-- A unified treatment, Proc. London Math. Soc. (3) {\bf 50}
(1985), 447 -- 479.


\bibitem{CG1} A.L.S. Corner and R. G\"obel, Small almost free modules
 with prescribed topological endomorphism rings, Rendiconti Sem.
 Mat. Univ. Padova {\bf 109} (2003) 217 -- 234.

\bibitem{DG} M. Dugas and R. G\"obel, Automorphism groups of fields II,
Commun. in Algebra, {\bf 25} (1997) 3777 -- 3785.

\bibitem{DG1} M. Dugas and R. G\"obel, Automorphism groups of geometric
lattices, Algebra Universalis {\bf 45} (2001), no. 4, 425--433.


\bibitem{EM} P. Eklof and A. Mekler, {\sl Almost free modules,
  set-theoretic methods} (revised edition), North-Holland, Elsevier,
  Amsterdam, 2002.

\bibitem{ES} P. Eklof and S. Shelah, Absolutely rigid systems and
  absolutely indecomposable groups, Abelian Groups and Modules,
  Trends in Math. (Birkh\"auser, 1999), 257 -- 268.

\bibitem{FrG} B. Franzen and R. G\"obel, The Brenner-Butler-Corner
theorem and its applications to modules, in: Abelian Group Theory
(Gordon \& Breach Science Publishers, London, 1987), 209 -- 227.

\bibitem{FK} E. Fried and J. Kollar, Automorphism groups of
fields, Colloqu. Math. Soc. Janos Bolyai {\bf 29} (1977), 293 --
317.

\bibitem{F} L. Fuchs, {\sl Abelian Groups}, vol. 2 (Academic
Press, 1973).

\bibitem{FG} L. Fuchs and R. G\"obel, Modules with absolute endomorphism
rings, submitted

\bibitem{G} R. G\"obel,  Vector spaces with five
distinguished subspaces, Results in Mathematics {\bf 11} (1987), 211
-- 228.

\bibitem{GM} R. G\"obel and W. May, Four submodules suffice for
  realizing algebras over commutative rings,  J. Pure Appl. Algebra
{\bf 65} (1990), 29 -- 43.

\bibitem{GM2} R. G\"obel and W. May, Independence in completions
and endomorphism algebras,  Forum Math.  {\bf 1} (1989), 215 -- 226.

\bibitem{GS} R. G\"obel and S. Shelah, Indecomposable almost free
modules|the local case, Canadian J. Math. {\bf 50} (1998), 719 --
738.


\bibitem{GT} R. G\"obel and J. Trlifaj, {\sl Endomorphism Algebras
and Approximations of Modules}, Walter de Gruyter Verlag, Berlin,
to appear (2006).

\bibitem{H} H. Heineken, Automorphism groups of torsionfree nilpotent
groups of class two.  Symposia Mathematica, {\bf 17} (1976), 235 --
250.


\bibitem{J} T. Jech, {\sl Set Theory}, Academic Press (1978).



\bibitem{Na} C. St. J. A. Nash-Williams, On well-quasi-ordering infinite trees,
Proc. Camb. Phil. Soc. {\bf 61} (1965), 697 -- 720.

\bibitem{S0} S. Shelah, Infinite abelian groups, Whitehead problem
and some constructions, Israel J. Math. {\bf 18} (1974), 243 --
256.

\bibitem{S} S. Shelah, Better quasi-orders for uncountable cardinals,
  Israel J. Math. 42 (1982), 177 -- 226.


\bibitem{Si} D. Simson, {\sl Linear Representations of Partially
Ordered Sets and Vector Space Categories},  Algebra,  Logic and
Applications, Vol. 4, Gordon \& Breach Science Publishers,
London, 1992.





\end{thebibliography}
\end{document}